\documentclass[12pts]{article}
\usepackage{times}
\usepackage{bm}
\usepackage{natbib}
\usepackage{amsmath}
\usepackage[figuresright]{rotating}

\usepackage{times}

\usepackage{amsfonts}
\usepackage{graphicx}

\newcommand{\N}{{{\cal N}}}

\newcommand{\Ob}{{{\cal O}}}

\newcommand{\Eet}{{\rm E}_{*}}
\newcommand{\Ee}{{\rm E}}
\newcommand{\var}{{\rm var}_{*}}

\newcommand{\dpg}[2]{\frac{\partial #1}{\partial #2}}
\newcommand{\ddpg}[2]{\frac{\partial^2 #1}{\partial #2^2}}

\newcommand{\si}{\sum_{i=1}^n}
\newcommand{\nn}{n^{-1}}
\newcommand{\tk}{\hat \theta}
\newcommand{\gtk}{g^{\hat \theta}}
\newcommand{\htk}{h^{\hat \gamma}}
\newcommand{\tki}{\hat \theta_{-i}}
\newcommand{\gtki}{g^{\hat \theta_{-i}}}
\newcommand{\UACV}{{{\rm UACV}}}
\newcommand{\AIC}{{ {\rm AIC}}}
\newcommand{\CRPS}{{{\rm CRPS}}}
\newcommand{\ECE}{{{\rm ECE}}}
\newcommand{\CV}{{{\rm CV}}}
\newcommand{\ER}{{{\rm ER}}}
\newcommand{\Fio}{\Phi_{\bar \Ob_n}}
\newcommand{\Fioi}{\Phi_{\bar \Ob_{n|i}}}
\newcommand{\Fioii}{\Phi_{\bar \Ob_{i}}}

\newcommand{\btau}{\mbox{\boldmath $\tau$}}
\newcommand{\bhl}{\mbox{\boldmath $\hat b$}(\theta)}

\newcommand{\bihl}{\hat b_i(\theta)}
\newcommand{\bb}{\mbox{\boldmath $b$}}
\newcommand{\hl}[1]{{\rm HL}(\theta,#1)}

\newcommand{\hlm}{{\rm hl}}

\newcommand{\phl}{{\rm PHL}_n(\theta)}

\def\Tr{\mbox{Trace}}

\newcommand{\KL}{{\rm KL}}
\newcommand{\EKL}{{\rm EKL}}

\newcommand{\tendD}{\longrightarrow \hspace{-0.50cm} ^D \hspace{0.50cm}}

\title{A universal approximate cross-validation criterion and its asymptotic distribution}

\author{Daniel Commenges $^{1,2}$ 
 \and C\'ecile Proust-Lima $^{1,2}$
\and C\'ecilia Samieri $^{1,2}$
\and Benoit Liquet $^{1,2}$}

\begin{document}





\maketitle
$^1$ INSERM, ISPED, Centre INSERM U-897-Epidemiologie-Biostatistique, Bordeaux,  F-33000\\
$^2$ Univ. Bordeaux, ISPED, Centre INSERM U-897-Epidemiologie-Biostatistique, Bordeaux,  F-33000, France\\ 
\vspace{5mm}

{\bf ABSTRACT.}
 A general framework is that the estimators of a distribution are obtained by minimizing a function (the estimating function) and they are assessed through another function (the assessment function). The estimating and assessment functions generally estimate risks. A classical case is that both functions estimate an information risk (specifically cross entropy); in that case Akaike information criterion (AIC) is relevant. In more general cases, the assessment risk can be estimated by leave-one-out crossvalidation. Since leave-one-out crossvalidation is computationally very demanding, an approximation formula can be very useful. A universal approximate crossvalidation criterion (UACV) for the leave-one-out crossvalidation is given.  This criterion can be adapted to different types of estimators, including penalized likelihood and maximum a posteriori estimators, and of assessment risk functions, including information risk functions and continuous rank probability score (CRPS). This formula reduces to Takeuchi information criterion (TIC) when cross entropy is the risk for both estimation and assessment. The asymptotic distribution of UACV and of a difference of UACV is given. UACV can be used for comparing estimators of the distributions of ordered categorical data derived from threshold models and models based on continuous approximations. A simulation study and an analysis of real  psychometric data are presented.

\vspace{2mm}


\newpage



\section{Introduction}

  Akaike information criterion (AIC) \citep{Aka73}  is useful for comparing parametric models.  AIC assumes parametric models and maximum likelihood estimators; it has been developed from the Kullback-Leibler divergence.
    \cite{Com08}  showed that a normalized difference of AIC of two models could be considered as estimating a difference of Kullback-Leibler risks of maximum likelihood estimators of the density based on the two models. Likelihood crossvalidation (LCV) has also been widely used for comparing parametric models. \cite{STO74} showed that LCV was asymptotically identical to AIC. LCV however is more flexible in that it can be applied to other estimators than maximum likelihood estimators (MLE), for instance to penalized likelihood estimators: see \cite{golub1979generalized,wahba1985comparison}. It has asymptotic optimality properties \citep{van2004asymptotic}. The leave-one-out crossvalidation is the most natural and one of the most efficient but it is also the most computationally demanding so that approximation formulas have been derived. \cite{STO74} was the first to give an approximation formula. Other works developed generalized  approximation crossvalidation for smoothing parameter selection for penalized splines \citep{xiang1996generalized,gu2001cross} or for penalized likelihood \citep{o1986statistical,Com07}. Cross-validation can also be applied to other assessment risks than Kullback-Leibler risk; for a review see \cite{Arl_Cel:2010:surveyCV}.

    We consider the following framework: estimators of the true density function are defined as minimizing an estimating function; the estimating function itself can be viewed as an estimator of a risk, that we call an ''estimating risk"; the estimators of the true density are assessed using an ''assessment risk". The assessment risk can be estimated by crossvalidation, allowing a choice between them. Leave-one-out crossvalidation is one of the best crossvalidation procedure but is very computationally demanding. The aim of this paper is to find a universal approximation for leave-one-out crossvalidation, valid whatever the estimating and assessment risks.

 Section 2 presents the framework and the universal approximate criterion (UACV); section 3 shows how it specializes to particular cases. In section 4 the asymptotic distributions of UACV  and of a difference of UACV are given. Section 5 presents a simulation study, while section 6 presents an illustration of the use of UACV for comparing estimators derived from threshold models and estimators obtained by  continuous approximations in the case of ordered categorical data with repeated measurements. Section 7 concludes.

\section{The universal approximate crossvalidation criterion}\label{derUACV}
\subsection{The risk and its estimation by crossvalidation}
Suppose that a sample of independently identically distributed (iid) variables $\bar \Ob_n=(Y_i, i=1,\ldots,n)$ is available. Based on $\bar \Ob_n$, an estimator $\gtk$ (where $\tk$ is short for $\hat \theta_n$) of the probability density function $f^*$ of the true distribution can be chosen in a family of distributions $(g^{\theta})_{\theta \in \Theta}$, $\Theta \subset \Re^p$. Estimators are chosen as minimizing an estimating function $\Phi_{\bar \Ob_n}(\theta)$ where $\Phi_{\bar \Ob_n}(\theta)=\nn \si \phi(\theta,Y_i)$: that is, $\tk= {\rm argmin} _{\theta} \Phi_{\bar \Ob_n}(\theta)$. We shall assume that $\phi(\theta,y)$ is a thrice differentiable function of $\theta$ for any $y$. Under mild conditions $\Phi_{\bar \Ob_n}$ converges toward $\Phi_{\infty}(\theta)=\Eet\{\phi(\theta,Y_i)\}$, where $\Eet$ means that the expectation is taken under the true probability specified by $f^*$. Thus, we can think of these estimating functions as estimators of a risk function, which by definition is  the expectation of a loss function. For instance if the loss is $\phi(\theta,Y)= -\log g^{\theta}(Y)$, $\Phi_{\bar \Ob_n}(\theta)$ estimates the risk $\Eet \{-\log g^{\theta}(Y)\}$; this risk is the cross entropy of $g^{\theta}$ relative to $f^*$ (see section \ref{AIC}). Maximum likelihood estimators (MLE) are derived this way \citep{Aka73,Com09}. A more general class of estimators of this form is that of M-estimators. Under some conditions given in \cite{van2000asymptotic}, $\tk$ converges in probability toward $\theta_0={\rm argmin} _\theta \Phi_{\infty}(\theta)$.
 Several estimating loss functions depend on $\theta$ only through $g^{\theta}$ but some may depend on $\theta$ directly, for instance through penalty terms.

Let us consider another loss function $\psi(g^{\theta},Y)$ which defines a risk function $R^{\psi}(g^{\theta})=\Eet\{\psi(g^{\theta},Y)\}$; we shall assume that  $\psi(g^{\theta},y)$ is a twice differentiable function in $\theta$ for all values of $y$. This risk function may serve for assessing estimators, but since estimators are random, the assessment risk is an expected risk: $\ER^{\psi}(\gtk)=\Eet\{\psi(\gtk,Y)\}$. Estimators with small assessment risks are preferred. The problem is to estimate the assessment risk (without knowing the true density $f^*$). This will be easier if  the loss itself does not involve $f^*$; this is the case for instance of $\psi(\theta,Y_i)=-\log g^{\theta}(Y_i)$ but not of the integrated squared error ${\rm ISE}(g^{\theta})=\int \{g^{\theta}(u)-f^*(u)\}^2du$. The expectation however is with respect to the true distribution, so the risk cannot be computed exactly but has to be estimated. If another sample $\bar \Ob'_n=(Y'_i, i=1,\ldots,n)$ iid with respect to $\bar \Ob_n$ were available, a natural estimator of the risk would be $\tilde \ER^{\psi}(\gtk)=\nn \si \psi(\gtk,Y'_i)$. This is an unbiased estimator of the risk.
This is often used by practitioners who split their original sample in a training and a validation sample. However this practice leads to a loss of efficiency since only half of the data are used for computing the estimator $\gtk$ and half of the data also for estimating its assessment risk. Crossvalidation estimators make a more efficient use of the information. In particular the leave-one-out crossvalidation criterion is:
$$\CV(\gtk)=\nn \si \psi(\gtki,Y_i),$$
where $\tki= {\rm argmin}~ \Phi_{\bar \Ob_{n|i}}$ and $\Phi_{\bar \Ob_{n|i}}=\frac{1}{n-1} \sum_{j\ne i}^n \phi(\theta,Y_j)$.
$\CV$ does nearly as well as if another sample $\bar \Ob'_n$ were available. Indeed it can immediately be seen that $\Ee\{\CV(\gtk)\}=\ER^{\psi}(g^{\hat \theta_{n-1 }})$. We shall see its asymptotic distribution in section \ref{s-asympt}. For comparing two estimators the difference of risks is relevant. This can of course be estimated by the difference of the estimated risks whose asymptotic distribution will also be given in section \ref{s-asympt}.

\subsection{The universal approximate crossvalidation criterion}
The leave-one-out crossvalidation criterion may be computationally demanding since it is necessary to run the maximization algorithm $n$ times for finding the $\tki, i=1,\ldots,n$ . For this reason an approximate formula is very useful. This is possible because $\tki$ is close to $\tk$. Indeed, under mild conditions given in \cite{van2000asymptotic}, $\sqrt{n} (\tk- \theta_0)$ has an asymptotic normal distribution; since this is also true for $\tki$,  this implies that $\tki-\tk=O_p(n^{-1/2})$. A Taylor expansion of $\dpg{\Fioi}{\theta}|_{\tki}$ around $\tk$ yields:
\begin{equation*}\label{1stepNR} \tki-\tk=-H^{-1}_{\Fioi}\dpg {\Fioi}{\theta}|_{\tk}+R_n, \end{equation*}
where $H_{\Fioi}=\ddpg{\Fioi}{\theta }|_{\tk}$, and $R_n$ is a quadratic form of $\tki-\tk$ involving second and third derivatives of $\Fioi$ taken in $\tilde \theta$ so that $||\tilde \theta_n-\tk|| \le ||\tki-\tk||$. Thus $||\tilde \theta_n-\tk||$ is also an $O_p(n^{-1/2})$. In virtue of the strong law of large numbers $\Fioi$ and its derivatives converge almost surely toward their expectations taken at $\theta_0$.
 By definition of $\Fio (\theta)$ we have the relation:
 \begin{equation} \label{relation} n\Fio (\theta)=(n-1)\Phi_{\bar \Ob_{n|i}}(\theta)+\phi(\theta,Y_i).\end{equation}
 Taking derivatives of the terms of this equation  and taking the values at $\tk$ we find $0=(n-1)\dpg {\Fioi}{\theta}|_{\tk}+\dpg{\phi(\theta,Y_i)}{\theta}|_{\tk}$ and we obtain that $\dpg {\Fioi}{\theta}|_{\tk}=-\hat d_i$ where
 \begin{equation}\hat d_i=\frac{1}{n-1}\dpg{\phi(\theta,Y_i)}{\theta}|_{\tk}.\end{equation}
 Note that on the right hand of the equation we could use a multiplicative factor  $\frac{1}{n}$ instead of $\frac{1}{n-1}$ since the difference between the two is a $O(n^{-2})$.
  Hence we have:
\begin{equation}\label{approxtheta1} \tki-\tk=H^{-1}_{\Fioi}\hat d_i +R_n, \end{equation}

Note that this implies that $\tki-\tk=O_p(\nn)$  because $H_{\Fioi}=O_p(1)$ and both $\hat d_i $ and $R_n$ are $O_p(n^{-1})$. But this in turn implies that $R_n$ is in fact an $O_p(n^{-2})$.
By twice derivating  (\ref{relation}) we obtain: $H_{\Fio} =\frac{n-1}{n}H_{\Fioi}+\frac{1}{n}H_{\Fioii}$ where $H_{\Fioii}=\ddpg{\Fioii}{\theta}|_{\tk}$; since the last term is an $O_p(n^{-1})$, we can replace $H_{\Fioi}$ by  $H_{\Fio}=\ddpg{\Fio}{\theta}|_{\tk}$ in (\ref{approxtheta1}) and obtain:
\begin{equation}\label{approxtheta} \tki-\tk=H^{-1}_{\Fio}\hat d_i + O_p(n^{-2}). \end{equation}

Developing now the assessment loss function for $\tki$ around $\tk$ yields:
\begin{equation*}\label{TaylorLoss} \psi(\gtki,Y_i)=\psi(\gtk,Y_i)+ (\tki-\tk)^T\hat v_i+O_p(n^{-2}), \end{equation*}
where
\begin{equation}\hat v_i=\dpg {\psi(g^{\theta},Y_i)}{\theta}|_{\tk}.\end{equation}
Replacing in this equation $\tki-\tk$ by its approximation in (\ref{approxtheta}) we obtain: $\psi(\gtki,Y_i)=\psi(\gtk,Y_i)+ \hat d_i^TH^{-1}_{\Fio}\hat v_i+O_p(n^{-2})$. Taking the mean of the left terms of these equations yields $\CV(\gtk)$ so that we have:

\begin{equation} \CV(\gtk)=\Psi(\gtk)+\Tr(H^{-1}_{\Fio}K)+O_p(n^{-2}),\end{equation}
where $\Psi(\gtk)=\nn \si \psi(\gtk,Y_i)$ and $K=\nn \si \hat v_i \hat d_i^T$.
We define UACV as:
\begin{equation}\label{UACV} \UACV(\gtk)=\Psi(\gtk)+\Tr(H^{-1}_{\Fio}K).\end{equation}

\subsection{Extension to incompletely observed data}
More generally the available observation is $\bar \Ob_n=(\Ob_i, i=1,\ldots,n)$ where $\Ob_i$ are sigma-fields containing observation of events generated by $Y_i$. Classical cases are that of right-censored observations common in survival analysis, left-censored observations common in biological measurements with detection limits, interval-censored data common in observations from cohort studies. Then risks functions of the type $\phi(\theta,Y_i)$ and $\psi(\gtk,Y_i)$ cannot be computed from observations (they are not $\bar \Ob_n$-measurable). The theory above does not apply and attempts to estimate such functions require that the model is well-specified, a very strong assumption in this context. A solution is to consider loss functions of the form $\phi(\theta,\Ob_i)$ and $\psi(\gtk,\Ob_i)$, to which the above results apply. The most natural choice is the log-likelihood. Such an approach is in fact heuristically taken when using AIC with censored data. Some theory and simulations have been given in \citet{Com07}, \citet{liquetchoice} and \citet{Com12prognosis} .

\section{Particular cases of UACV}
\subsection{Maximum likelihood estimators and information risk: AIC} \label{AIC}
Suppose we take: $\phi(\theta,Y_i)=\psi(g^{\theta},Y_i)=-\log g^{\theta}(Y_i)$. Then, the estimating function is minus the loglikelihood which estimates the estimating risk, here the cross-entropy \citep{cov91} of $g^{\theta}$ with respect to the true density $f^*$: $\Eet\{- \log g^{\theta}(Y)\}=H(f^*)+\KL(g^{\theta};f^*)$, where $H(f^*)=-\Eet\{\log f^*(Y)\}$ is the entropy of $f^*$ and $\KL(g^{\theta};f^*)=\Eet\{\log \frac{f^*(Y)}{g^{\theta}(Y)}\}$ the Kullback-Leibler divergence of $g^{\theta}$ relative to $f^*$.
 As for the assessment risk, this is the expected cross entropy:
 \begin{equation} \label{ECE} \ECE(g^{\tk})=\Eet[\Eet\{- \log \gtk (Y)|\bar \Ob_n\}]=H(f^*)+\EKL(g^{\tk};f^*),\end{equation}
  where $\EKL(g^{\tk};f^*)=\Eet\{\log \frac{f^*(Y)}{\gtk(Y)}\}$. In that case the loss functions for estimating and assessment are the same. The estimating and assessment risks are nearly the same; there is however a dissymmetry in that the estimating risk is a cross entropy while, because $\gtk$ is random, the assessment risk is an expected cross-entropy.

 In that case the leading term of (\ref{UACV}) is minus the maximized (normalized) loglikelihood. We have also that $\hat v_i$ is the individual score and $\hat d_i=\frac{1}{n-1}\hat v_i$ so that UACV is identical to a normalized version of Takeuchi information criterion (TIC). If the model is well specified  $K$ tends in probability toward $I_L$, where $I_L$ is the individual information matrix. The Hessian $H^{-1}_{\Fio}$ also tends toward $I_L$ so that the correction term tends toward $p$, the number of parameters. Thus, if the model is not too badly specified, TIC is approximately equal to AIC. We have $\UACV=\frac{1}{2n} TIC\approx \frac{1}{2n} AIC$, and this estimates the {\it expected} cross-entropy of the estimator.

\subsection{Maximum a posteriori and maximum penalized likelihood estimators}
Estimators can be obtained by maximizing the maximum a posteriori (MAP) density of $\theta$: $\si \log g^{\theta}(Y_i)-J(\theta)$, where $J(\theta)$ is minus the log of the prior density of $\theta$. This corresponds to a loss function $\phi(\theta,Y)=-\log g^{\theta}(Y)+\nn J(\theta)$. Here the loss function depends on $n$ but since the MAP estimators also converge (toward the same limit as the MLE), the above results still hold so that UACV can be applied to the MAP estimators. The assessment risk can be the expected cross entropy. Thus, in that case the estimating risk and the assessment risk are clearly different.

In the case where at least one of the parameters is a function, penalized likelihood can be applied for obtaining a smooth estimator \citep{o1986statistical}. The function is approximated on a spline basis so that the problem reduces to a parametric one with a large number of parameters. This leads to $\phi(\theta,Y_i)=-\log g(Y_i;\theta) + J(\theta,n)$, where $J(\theta,n)$ is a penalty, which involves for instance the norm of the second derivative of the function and which may depend on $n$. $J(\theta,n)$ may depend on $n$ in such a way that the estimator converges. The assessment risk can be the expected cross entropy as in section \ref{AIC}. In that case the leading term in (\ref{UACV}) is minus the (normalized) loglikelihood taken at the value of the penalized likelihood estimator, $\hat \theta$. We still have that $\hat v_i$ is the individual score (although taken at the value of the penalized likelihood estimator), and $\hat d_i=\frac{1}{n-1}(\hat v_i+\dpg{J(\theta,n)}{\theta}|_{\hat \theta})$. We have that $H_{\Fio}= -\ddpg {L_{\bar \Ob_n}}{\theta}_{|\tk}+\ddpg{J(\theta,n)}{\theta}_{|\tk}$. In the case where $\ddpg{J(\theta,n)}{\theta}$ is a definite positive matrix, we have $H_{\Fio}> -\ddpg {L_{\bar \Ob_n}}{\theta}_{|\tk}$ and $H^{-1}_{\Fio}< -\left (\ddpg {L_{\bar \Ob_n}}{\theta}_{|\tk}\right )^{-1}$ so that the correction term is smaller than for MLE. The correction term is often interpreted as the equivalent number of parameters \citep{Com07}.

\subsection{Hierarchical likelihood estimators}
Let us consider the following model: conditionally on $b_i$, $Y_i$ has a density $g_{Y|b}(.;\theta, b_i)$, where $b_i$ are random effects (or parameters). The $(Y_i,b_i)$ are iid, where $Y_i$ is multivariate of dimension $n_i$. The $b_i$ are assumed to have density $g_b(.; \btau)$ with zero expectation, where $\btau$ is a vector of parameters. Typically $Y_i$ is (at least partially) observed while $b_i$ is not.

 Estimators of both $\theta$ and $\bb=(b_1,\ldots,b_n)$ are defined as maximizing the h-loglikelihood which can be written:
$\hl{\bb,\btau}=\frac{1}{n}\sum _{i=1}^n \hlm(Y_i;\theta,b_i, \btau)$ with $\hlm(Y_i;\theta,b_i,\btau)=g_Y(Y_i;\theta, b_i)+\log g_b(b_i; \btau)$.

\cite{commenges2010inference} showed (via a profile likelihood argument) that the hierarchical likelihood estimators for $\theta$ are M-estimators: $\phl=\hl{\bhl}=\frac{1}{n}\sum _{i=1}^n \hlm(Y_i;\theta,\bihl)$. Thus the estimating loss is $\hlm(Y_i;\theta,\bihl)$. Estimators may be compared using UACV. If we use the expected cross entropy for assessing the estimators, it is necessary to compute the marginal likelihood in $\hat \theta$ and the values of the individual scores. Of course, this is precisely what we attempted to avoid by using hierarchical likelihood; we have however to do it only once (rather than repeatedly within an iterative algorithm), so it may be computationally feasible. For the correcting term we have also to compute the Hessian of the hierarchical likelihood and the individual scores.

\subsection{Restricted AIC}
\cite{liquetchoice} have proposed a modification of AIC and LCV for assessing estimators based on information $\bar \Ob_n$ while the assessment risk is based on a smaller information  $\Ob'_{n+1} \subset \Ob_{n+1}$. More specifically, the estimator is based on the sample $\bar \Ob_n=(Y_i, i=1,\ldots,n)$ but the assessment risk is based on a random variable $Z$ which is a coarsened version of $Y$. For instance $Z$ is a dichotomization of $Y$: $Z=1_{Y>l}$.
For this case, the restricted AIC (RAIC) was derived by both direct approximation of the risk and by approximation of the LCV. RAIC is clearly a particular case of $\UACV$ for the case: $\phi(\theta,Y_i)=-\log g^{\theta}(Y_i)$ and  $\psi(g^{\theta}(Y_i))=-\log g^{\theta}(Z_i)$.

\subsection{Case where the estimators are assessed by CRPS}
\cite{gneiting2007strictly} studied scoring rules and particularly the continuous rank probability score (CRPS). Its inverse can be used as a loss function and is defined as:
$$\CRPS^*(G(.,\theta),Y)= \int_{-\infty}^{+\infty} \{G(u,\theta)-1_{u\ge Y}\}^2~du,$$
where $G(.,\theta)$ is the cumulative distribution function (cdf) of a distribution in the model. The risk is a Cramer-von Mises type distance : $d(G,G^*)=\int \{G(u)-G^*(u)\}^2~du$.
It may be interesting to assess MLE's using this assessment risk.
In $\UACV$, the leading term is $\nn \si \CRPS^*(G(.,\tk),Y_i)$; for the correcting term, $H_{\Fio}$ is the Hessian of the loglikelihood and $K$ must be computed with $\hat v_i=\dpg {\psi}{\theta}|_{\tk}=2\int_{-\infty}^{+\infty} \{G(u,\tk)-1_{u\ge Y}\}\dpg{G(u,\theta)}{\theta}|_{\tk}~du.$

\subsection{Estimators based on continuous approximation of categorical data}\label{categorical}
Assume $Y$ is an ordered categorical variable taking values $l=0,1,\ldots,L$. Here for simplicity we consider that $Y$ is univariate. Several models are available for this type of variables.
Threshold link models assume that $Y_i=l$ if a latent variable $\Lambda_i$ takes values in the interval $(c_l,c_{l+1})$ for $l=1,\ldots,L$, with $c_{L+1}=+\infty$ and $Y_i=0$ if $\Lambda_i < c_1$:
\begin{equation}\label{ThresholdLink} Y_i=\sum_{l=1}^{L+1} 1_{\{\Lambda_i\in(c_l,c_{l+1})\}}l.\end{equation}
 $\Lambda_i$ itself can be modeled as a noisy linear form of explanatory variables $\Lambda_i=\beta x_i+ \varepsilon_i$, with $\varepsilon_i$ having a normal distribution of mean zero and variance $\sigma^2$, and where $x_i$ are explanatory variables. The parameters are $\theta=(c_1,\ldots,c_L, \beta, \sigma)$. An estimator of the distribution can be obtained by maximum likelihood leading to define $g^{\hat \theta}$. For assessing the estimator it is natural to use $\ECE(g^{\hat \theta})$ (defined in \ref{ECE}). Note that since $Y$ is discrete the densities are defined with respect to a counting measure that is,  $g^{\hat \theta}(l)$  defines the probability that $Y=l$.

One may also make a continuous approximation which leads to simpler computations and may be more parsimonious, especially if $Y$ is multivariate as in the illustration of section \ref{illustration}. In this approach we consider the model $Y_i= \beta x_i+ \varepsilon_{i}$.  Maximizing the likelihood of this model for observations of $Y_i$ leads to a probability measure specified by the density $h_c^{\hat \gamma}$. This is however a density relative to Lebesgue measure. This probability measure gives zero probabilities to $\{Y_i=l\}$ for all $l$, and this yields infinite value for ECE (meaning strong rejection of this estimator). However from $h_c$ a natural estimator of $f^*$ can be constructed by gathering at $l$ the mass around $l$: $h^{\hat \gamma}(l)=\int_{l-1/2}^{l+1/2} h_c^{\hat \gamma}(u)~du$, for $l=1,\ldots,L-1$, and $h^{\hat \gamma}(0)=\int_{-\infty}^{1/2} h_c^{\hat \gamma}(u)~du$, $h^{\hat \gamma}(L)=\int_{L-1/2}^{+\infty} h_c^{\hat \gamma}(u)~du$. UACV can be computed for this estimator for estimating its ECE. The first term of $\UACV(h^{\hat \gamma})$ can be interpreted as the loglikelihood obtained by this estimator with respect to the counting measure. For the correcting term we need the Hessian of the loglikelihood of $h_c^{\hat \gamma}$ and we have to compute $\hat v_i=\dpg {\psi(h^{\gamma},Y_i)}{\gamma}|_{\tk}$. For instance if $Y_i=l$ for $l=1\ldots,L-1$ we have
$$\hat v_i=-\frac{\int_{l-1/2}^{l+1/2} \dpg{h_c^{\hat \gamma}}{\gamma}(u)~du}{\int_{l-1/2}^{l+1/2}h_c^{\hat \gamma}(u)~du}.$$
Since the denominator is the probability under $h_c^{\hat \gamma}$ that  $Y\in (l-1/2,l+1/2)$, $\hat v_i$ can be interpreted as the conditional expectation (under $h_c^{\hat \gamma}$) of the individual score. Thus if $h_c^{\hat \gamma}$ does not vary much on $(l-1/2,l+1/2)$, $\hat v_i$ is close to $-(n-1)\hat d_i$. Using the same arguments as in section \ref{AIC} we obtain that UACV is close to correcting by the number of parameters as in AIC, a criterion that we call AIC$_d$. This is what \cite{proust2011misuse} proposed, and this is likely to be a good approximation if the number of modalities of $Y$ is large.

\section{Asymptotic distribution and tracking interval}\label{s-asympt}
\cite{Com08} using results of \cite{vuong1989likelihood} studied the asymptotic distribution of a normalized difference of AIC as an estimator of a difference of Kullback-Leibler risks. Here similar arguments are applied to study the asymptotic distribution of UACV and a difference of UACV. Since UACV and CV differ by an $O_p(n^{-2})$, this will also give the asymptotic distribution of CV and a difference of CV. By the continuous mapping theorem, the asymptotic distribution of $\UACV(\gtk)$ is the same as that of $\Psi(g^{\theta_0})$. Since the latter quantity is a mean, it immediately follows by the central limit theorem that:
\begin{equation}n^{1/2}\{\UACV(\gtk)-R^{\psi}(g^{\theta_0})\} \tendD \N (0,\kappa_*^2),\end{equation}
where $\kappa_*^2= \var \psi(g^{\theta_0},Y)$. We can also write:
\begin{equation}n^{1/2}\{\UACV(\gtk)-\ER^{\psi}(\gtk)\} \tendD \N (0,\kappa_*^2),\end{equation}
and $\kappa_*^2$ can be estimated by the empirical variance of $\psi(\gtk,Y_i)$, $i=1,\ldots,n$.

 Given two estimators, $\gtk$ and $\htk$, it is interesting to estimate the difference of their risks: $\Delta^{\psi}(\gtk,\htk)=\ER^{\psi}(\gtk)-\ER^{\psi}(\htk)$. The obvious estimator is: $D_{\UACV}(\gtk,\htk)=\UACV(\gtk)-\UACV(\htk).$
We focus on the case where $g^{\beta_0} \ne h^{\gamma_0}$. We obtain in that case using the same arguments as above:
\begin{equation} \label{asympt} n^{1/2}\{D_{\UACV}(g^{\hat \beta_n},h^{\hat \gamma_n})-\Delta(g^{\hat \beta_n},h^{\hat \gamma_n})\} \tendD \N (0,\omega_*^2),\end{equation}
where $\omega_*^2={\rm var} \left \{ \psi(g^{\theta_0},Y)-\psi(h^{\gamma_0},Y)\right\}$, and this can be estimated by the empirical variance of $\left \{ \psi(\gtk,Y_i)-\psi(\htk,Y_i)\right\}.$

From this we can compute the so called tracking interval $(A_n,B_n)$, where $A_n=D_{\UACV}(g^{\hat \beta_n},h^{\hat \gamma_n})-z_{\alpha/2}n^{-1/2}\hat \omega_n$ and $B_n=D_{\UACV}(g^{\hat \beta_n},h^{\hat \gamma_n})+z_{\alpha/2}n^{-1/2}\hat \omega_n$, where $z_u$ is the $u^{th}$ quantile of the standard normal variable.

Note that $\omega_*$ is in general much lower than $\kappa_*$. This has been shown by \citet{Com07} for the cross entropy risk and comes from the fact that $\psi(\gtk,Y_i)$ and $\psi(\htk,Y_i)$ are often positively correlated.

\section{Simulation: choice of estimators for ordered categorical data}

\subsection{Design}
We conducted a simulation study to illustrate the use of UACV for comparing estimators derived from threshold link models and estimators obtained by a continuous approximation in the case of ordered categorical data (see section 3.6). The
aim was to assess the performance of UACV as an estimator of ECE, and to compare it to the normalized
 naive AIC criterion (noted AIC) and the normalized  AIC criterion computed on the counting measure (noted AIC$_d$). Performances of these criteria were studied in the case where the number of modalities ($L+1$) of the response variable $Y$ is small (section \ref{section-small}), and where $L$ is large (section \ref{section-large}).

{\bf True distributions}\\
For all the simulations, the data came from a threshold model specified by:
\[
   \left\{
  \begin{array}{ll}
    Y_i=l &  \textrm{if $\Lambda_i\in [c_l;c_{l+1})$ for $l=1,\ldots,L$}\\
    Y_i=0 & \textrm{if $\Lambda_i<c_1$ and $c_{L+1}=+\infty$ }\\
   \multicolumn{2}{l}{\Lambda_i=\beta_0+\beta_1X^1_i+ \beta_2X^2_i+\varepsilon_i \ \ \ i=1,\ldots,n}\\
   \end{array} \right.
\]

where $\varepsilon_i$ had a normal distribution of mean zero and variance $\sigma^2$ and the two explanatory variables $X^1$ and $X^2$ came from a standard normal distribution. In order to not disadvantage the linear continuous approximation compared to the threshold model, the parameters $c_1,\ldots,c_{L}$ were chosen as the solution of the following equations:
\[
   \left\{
  \begin{array}{l }
    P(Y_i<c_1)=P(Y_i>c_L)\\
	P(Y_i<c_1)=P(c_1<Y_i<c_2),\\
	c_{i+1}=c_i+m \text{ with } m=(c_{L}-c_1)/(L-1) \\
 \end{array} \right.
\]

{\bf The different models}\\
For each generated sample, we fitted the threshold model as previously defined,
and a linear model assuming a linear continuous approximation of the response variable $Y$, $Y_i=\beta_0+\beta_1X^1_i+ \beta_2X^2_i+\varepsilon_i$. Both models were estimated by maximum likelihood.

For all simulations, N = 1000 samples were generated. The true criterion ECE (available by simulation) was computed by  Monte-carlo approach.

\subsection{Small number of modalities ($L=4$)}\label{section-small}
We consider here the case where the number of modalities of $Y$ is relatively small ($L+1=5$). In this simulation, we fixed $\beta_0=1$,
$\beta_1=-2.1$, $\beta_2= -3.7$ and $\sigma^2= 4$. In Table \ref{small.level}, we present, for different sample sizes $n$, the results for the different empirical criteria AIC, AIC$_d$ and UACV which can be compared to the  true criterion ECE. For any sample size, the threshold model provided a
better ECE than the linear model (positive difference). It appeared that
UACV had a small bias for all the sample sizes (of order $10^{-3}$). The two others criteria
AIC and AIC$_d$ were also in favor of a threshold model. However, as expected, the naive normalized AIC failed to estimate the true ECE risk due to the wrong probability measure (Lebesgue measure instead of a counting measure). We can note that the criterion AIC$_d$ estimated
ECE relatively well, with a bias around $10^{-2}$ and $10^{-3}$.

\begin{table}[ht]
\begin{center}
\caption{Performance of the criteria for a small number of modalities ($L+1=5$) and different sample sizes. Mean over 1000 replications of the
difference of criteria UACV, AIC$_d$, AIC. ECE is the true risk. Biais of the criteria as estimator of ECE} \label{small.level}
\begin{tabular}{lccccccc}
  \hline
 & ECE & UACV & AIC$_d$ & AIC & Bias UACV & Bias AIC$_d$ & Bias AIC \\
  \hline
\multicolumn{8}{l}{{\bf n=300}}\\
Linear & 1.2650 & 1.2661 & 1.2702 & 1.5998 & 0.0011 & 0.0053 & 0.3348 \\
  Threshold & 0.9874 & 0.9835 & 0.9836 & 0.9836 & -0.0039 & -0.0038 & -0.0038 \\
  Difference  & 0.2775 & 0.2826 & 0.2866 & 0.6162 & 0.0051 & 0.0091 & 0.3387 \\
\multicolumn{8}{l}{{\bf n=500}}\\
Linear & 1.2594 & 1.2635 & 1.2660 & 1.6468 & 0.0041 & 0.0066 & 0.3875 \\
  Threshold & 0.9753 & 0.9810 & 0.9810 & 0.9810 & 0.0057 & 0.0057 & 0.0057 \\
  Difference  & 0.2840 & 0.2825 & 0.2849 & 0.6658 & -0.0016 & 0.0009 & 0.3818 \\
 \multicolumn{8}{l}{{\bf n=3000}}\\
Linear & 1.2612 & 1.2617 & 1.2621 & 1.6325 & 0.0005 & 0.0009 & 0.3713 \\
  Threshold & 0.9763 & 0.9749 & 0.9749 & 0.9749 & -0.0014 & -0.0014 & -0.0014 \\
  Difference  & 0.2850 & 0.2868 & 0.2873 & 0.6577 & 0.0019 & 0.0023 & 0.3727 \\
   \hline
\end{tabular}
\end{center}
\end{table}

\subsection{Large number of modalities} \label{section-large}
We consider here the case where the number of modalities of $Y$ is  relatively large ($L+1=20$). In this simulation, we fixed $\beta_0=1$,
$\beta_1=-0.3$, $\beta_2= -1.7$ and $\sigma^2= 4$. The results of this simulation are presented in Table \ref{large.level}.
For any sample size, the linear model provided a
better ECE than the threshold model (negative difference). It appeared that
UACV had a small bias for all the sample sizes (of order $10^{-3}$ and $10^{-4}$). The AIC$_d$ criterion gave similar results as the UACV criterion while the AIC criterion failed to find the best estimator (positive difference).

\begin{table}[ht]
\begin{center}
\caption{Performance of the criteria for a large number of modalities ($L+1=20$) and different sample sizes. Mean over 1000 replications of the
difference of criteria UACV, AIC$_d$, AIC. ECE is the true risk. Biais of the criteria as estimator of ECE} \label{large.level}
\begin{tabular}{lccccccc}
  \hline
 & ECE & UACV & AIC$_d$ & AIC & Bias UACV & Bias AIC$_d$ & Bias AIC \\
  \hline
\multicolumn{8}{l}{{\bf n=300}}\\
Linear & 2.6794 & 2.6792 & 2.6796 & 2.8076 & -0.0002 & 0.0002 & 0.1282 \\
Threshold & 2.7104 & 2.7070 & 2.7069 & 2.7069 & -0.0035 & -0.0036 & -0.0036 \\
Difference   & -0.0311 & -0.0278 & -0.0273 & 0.1007 & 0.0033 & 0.0038 & 0.1317 \\
 \multicolumn{8}{l}{{\bf n=500}}\\
Linear& 2.6765 & 2.6740 & 2.6742 & 2.7117 & -0.0026 & -0.0024 & 0.0352 \\
Threshold& 2.6933 & 2.6899 & 2.6899 & 2.6899 & -0.0033 & -0.0034 & -0.0034 \\
Difference & -0.0167 & -0.0160 & -0.0157 & 0.0218 & 0.0008 & 0.0010 & 0.0385 \\
 \multicolumn{8}{l}{{\bf n=3000}}\\
Linear & 2.6736 & 2.6716 & 2.6717 & 2.7107 & -0.0020 & -0.0019 & 0.0371 \\
Threshold& 2.6746 & 2.6725 & 2.6725 & 2.6725 & -0.0021 & -0.0021 & -0.0021 \\
Difference& -0.0010 & -0.0009 & -0.0008 & 0.0382 & 0.0001 & 0.0001 & 0.0391 \\
   \hline
\end{tabular}
\end{center}
\end{table}

\section{Illustration on the choice of estimators for psychometric tests}\label{illustration}

In epidemiological studies, cognition is measured by psychometric tests which usually consist in the sum of items measuring one or several cognitive domains. A common example is the Mini Mental State Examination (MMSE) score \citep{folstein1975}, computed as the sum of 30 binary items (grouped in 20 independent items) evaluating memory, calculation, orientation in space and time, language, and word recognition; for this reason it is called a ''sumscore" and ranges from 0 to 30. Although in essence psychometric tests are ordered categorical data, they are most often analyzed as continuous data. Indeed, they usually have a large number of different levels and, especially in longitudinal studies, models for categorical data are numerically complex. Recently, \cite{proust2011misuse} defined a latent process mixed model to analyse repeated measures of discrete outcomes involving either a threshold model or an approximation of it using continuous parameterized increasing functions. Comparison of models assuming either categorical data (using the threshold model) or continuous data (using continuous functions) was done with an AIC$_d$, computed with respect to the counting measure. In this illustration, we use UACV to compare such latent process mixed models assuming either continuous or ordered categorical data when applied on the repeated measures of the MMSE and its calculation subscore in a large sample from a French prospective cohort study.

\subsection{Latent process mixed models}
In brief, the latent process mixed model assumes that a latent process $(\Lambda_i(t))_{t \geq 0}$ underlies the repeated measures of the observed variable $Y_i(t_{ij})$ for subject $i$ ($i=1,...,n$) and occasion $j$ ($j=1,...,n_i$). First, the latent process trajectory is defined by a standard linear mixed model \citep{laird1982}: $\Lambda_i(t)=X_{i}(t)^T \beta + Z_i(t)^Tb_i +\epsilon_i(t_{ij})$ for $t \geq 0$ where $X_{i}(t)$ and $Z_i(t)$ are distinct vectors of time-dependent covariates associated respectively with the vector of fixed effects $\beta$ and the vector of random effects $b_i$ ($b_i \sim MVN(\mu,D)$), and $\epsilon_i(t) \sim \mathcal{N}(0,\sigma^2)$. We further assume that $b_{i0}$, the first component of $b_i$ that usually represents the random intercept, is $N(0,1)$ for identifiability; except for the variance of  $b_{i0}$, $D$ is an unstructured variance matrix.

 Then, a measurement model links the latent process with the observed repeated measures.
 For ordered categorical data,  a standard threshold model as defined in (\ref{ThresholdLink}) (section \ref{categorical}) for the univariate case is well adapted.
For continuous data,  the link has been modeled as $H(Y_i(t_{ij});\eta)=\Lambda_i(t_{ij})$ where $H(.;\eta)$ is a monotonic increasing transformation. Three families of such transformations are considered: (i) $H(y;\eta) = \dfrac{h(y;\eta_1,\eta_2)-\eta_3}{\eta_4}$ where $h(.;\eta_1,\eta_2)$ is the Beta cdf with parameters $(\eta_1,\eta_2)$; (ii) $H(y;\eta) = \eta_1 + \sum_{l=2}^{m+2} \eta_{l} B^I_l(y)$ where $(B^I_l)_{l=2,m+2}$ is a basis of quadratic I-splines with $m$ nodes; (iii)  $H(y;\eta) = \dfrac{y-\eta_1}{\eta_2}$ which gives the standard linear mixed model.

Latent process mixed models are estimated within the Maximum Likelihood framework using \texttt{lcmm} function of \texttt{lcmm} R package. When assuming continuous data, the log-likelihood can be computed analytically using the Jacobian of $H$. In contrast, when assuming ordered categorical data, an integration over the random-effects distribution intervenes in the likelihood computation which is approximated by a Gauss-Hermite quadrature. Further details can be found in \cite{proust2006} and \cite{proust2012psycho}.

UACV is computed from the loglikelihood $\Psi$ obtained for the maximum likelihood estimators $\hat{\theta}$ with respect to the counting measure :

\begin{equation}\label{Loglik_appli}
\begin{split}
\Psi(\hat{\theta}) &= \sum_{i=1}^n \int_{b_i} \prod_{j=1}^{n_i} P(Y_{ij}|b_{i}) f_b(b_i)db_i \\
&= \sum_{i=1}^n \int_{b_i} \prod_{j=1}^{n_i} \prod_{l=0}^{L} \left ( P(Y_{ij}=l |b_{i}) \right )^{I_{Y_{ij}=l}} f_b(b_i)db_i \\
&= \sum_{i=1}^n \int_{b_i} \prod_{j=1}^{n_i} \prod_{l=0}^{L} \left ( P(c_{l} \leq \Lambda_i(t_{ij})+\epsilon_i(t_{ij}) < c_{l+1} |b_{i}) \right )^{I_{Y_{ij}=l}}  f_b(b_i)db_i
\end{split}
\end{equation}

where $\Lambda_i(t_{ij})+\epsilon_i(t_{ij}) \sim N(X_{i}(t)^T \beta + Z_i(t)^T\mu,Z_i(t)^TD Z_i(t)+\sigma^2)$, $c_0=-\infty$, $c_{L+1}=+\infty$, and either $c_l$ ($l=1,...,L$) are the estimated thresholds when a threshold model is considered, or $c_l=H(l-\dfrac{1}{2},\hat{\eta})$ ($l=1,...,L$) when monotonic increasing families of transformations are used.

\subsection{Application: categorical psychometric tests}

Data come from the French prospective cohort study PAQUID initiated in 1988 to
study normal and pathological aging (\cite{letenneur1994}). Subjects included
in the cohort were 65 and older at initial visit and were followed up to 10
times with a visit at 1, 3, 5, 8, 10, 13, 15, 17 and 20 years after the initial visit. At
each visit, a battery of psychometric tests including the MMSE was completed. In the present analysis, all the subjects free of dementia at the 1-year visit and who had at least one MMSE measure during the whole follow-up were included. Data from baseline were removed to avoid modeling the first-passing effect. A total of 2914 subjects were included with a median of 2 (Interquartile range=3-5) repeated measures of MMSE sumscore. The distributions of MMSE sumscore ranging from 0 to 30, and of its calculation subscore, ranging from 0 to 5, are displayed in Figure \ref{histo_MMSE_CALC}.

\begin{figure}
\centering
\includegraphics[width=1.0\textwidth]{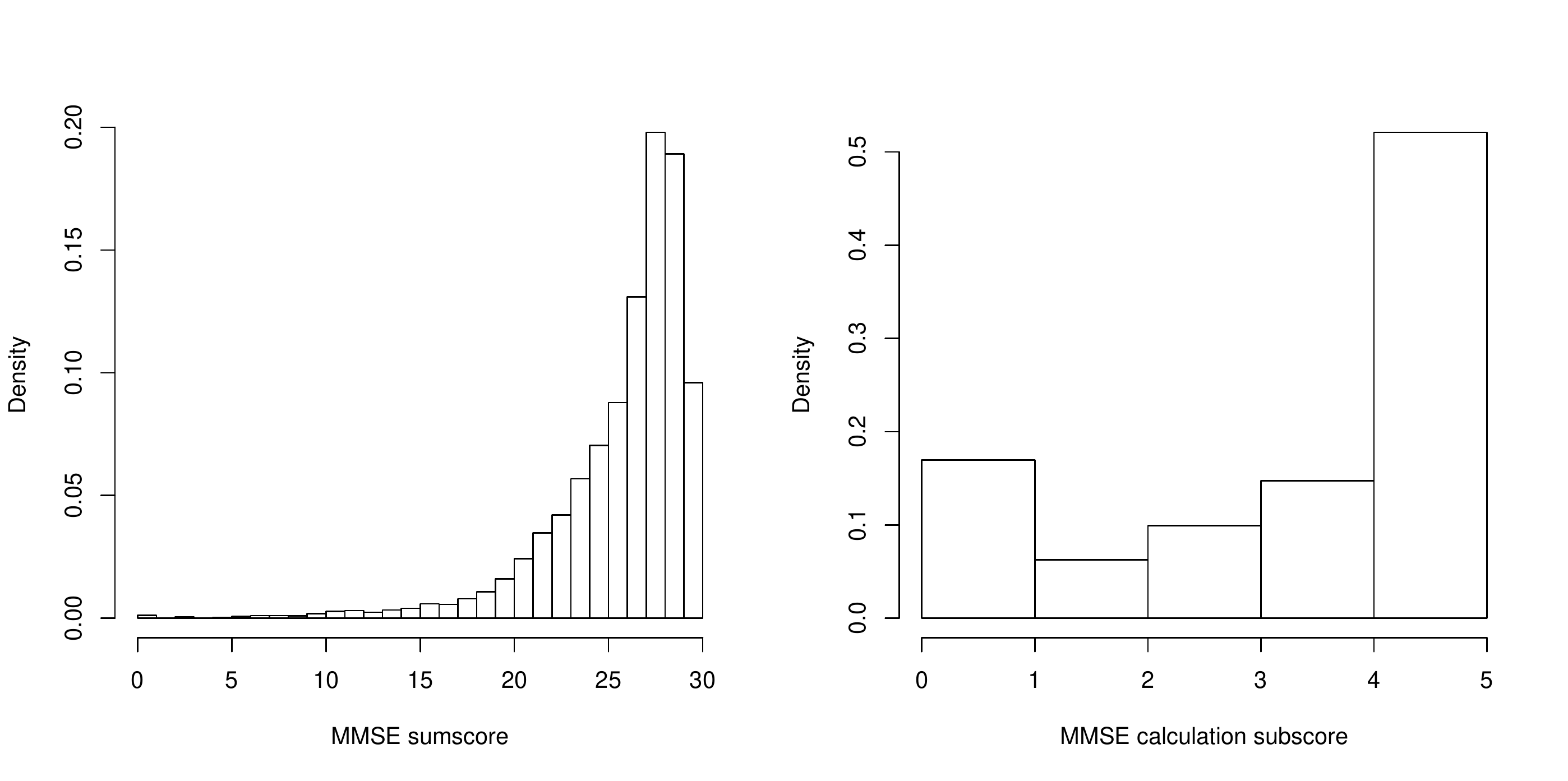}
\caption{Distributions of MMSE Sumscore and MMSE Calculation Subscore in the PAQUID sample (n=2,914). Data were pooled from all available visits for a total of 10,846 observations.}
\label{histo_MMSE_CALC}
\end{figure}

The trajectory of the latent process was modeled as an individual quadratic function of age with correlated random intercept, slope and quadratic slope ($Z_i(t)=(1,\text{age}_i,\text{age}_i^2)$), and an adjustment for binary covariates educational level (EL=1 if the subject graduated from primary school) and gender (SEX=1 if the subject is a man) plus their interactions with age and quadratic age (so that $X_i(t)=Z_i(t)\otimes(1,\text{EL}_i,\text{SEX}_i)$). For MMSE sumscore, in addition to the threshold model, the linear, Beta cdf and I-splines (with 5 equidistant nodes) continuous link functions were considered. For calculation subscore, in addition to the threshold model, only the linear link was considered.

 \subsection{Results}
  Table \ref{table_UACV} gives the assessment citeria for estimators based on the different models, and table \ref{table_diffUACV} provides the differences in UACV or AIC and their 95\% tracking interval. For the MMSE sumscore, the mixed model assuming the standard linear transformation yielded a clearly worse UACV than other models accounting for nonlinear relationships with the underlying latent process. The model involving a Beta cdf gave a similar risk as the one involving the less parsimonious I-splines transformation ($D_{\UACV}{\rm (Beta ~~cdf, I-splines)}= -0.0070$, 95\% Tracking interval:  $[-0.0152,0.0012]$). Finally, the mixed model considering a threshold link model, which is numerically demanding, gave the best fit but remained relatively close to the simpler ones assuming a Beta cdf ($D_{\UACV} {\rm (Beta ~~cdf,Thresholds)} = 0.0200$, 95\% Tracking interval: $[0.0097,0.0303]$) or a I-splines transformation ($D_{\UACV}  {\rm (I-splines, Thresholds)}= 0.0270$, 95\% Tracking interval: $[0.0166,0.0374]$). For the interpretation of these values \cite{Com08} suggested to qualify values of order $10^{-1}$, $10^{-2}$ and $10^{-3}$ as ''large", ''moderate" and ''small" respectively; moreover for multivariate observations, it was suggested to divide by the total number of observations rather by the number of independent observations. With this correction (which amounts to divide the current values by a factor of $3.7=10846/2914$) the differences between the linear model and the other models can be qualified as ''large", and the differences between the threshold model and both beta cdf and I-splines are between ''moderate" and ''small".  Figure 2  displays the estimated link functions in (A) and the predicted mean trajectories of the latent process according to educational level in (B) from the models involving either a linear, a beta cdf, a I-splines or a threshold link function. The estimated link functions as well as the predicted trajectories of the latent process are very close when assuming either  beta cdf,  I-splines or a threshold link function but they greatly differ when assuming a linear link. This difference also shows up in the effects of covariates with associations distorted when not accounting for nonlinear transformations (as demonstrated in \cite{proust2011misuse}). For example in this application, a significant overall increase with age of the difference between educational levels is found when assuming a linear model (p=0.011 for interaction with age and age squared) while it is a significant overall decrease with age which is found in the other models (p-value $< 0.001$ for interaction with age and age squared).

\begin{figure}\label{estimatedlink}
\centering
\includegraphics[width=1.0\textwidth]{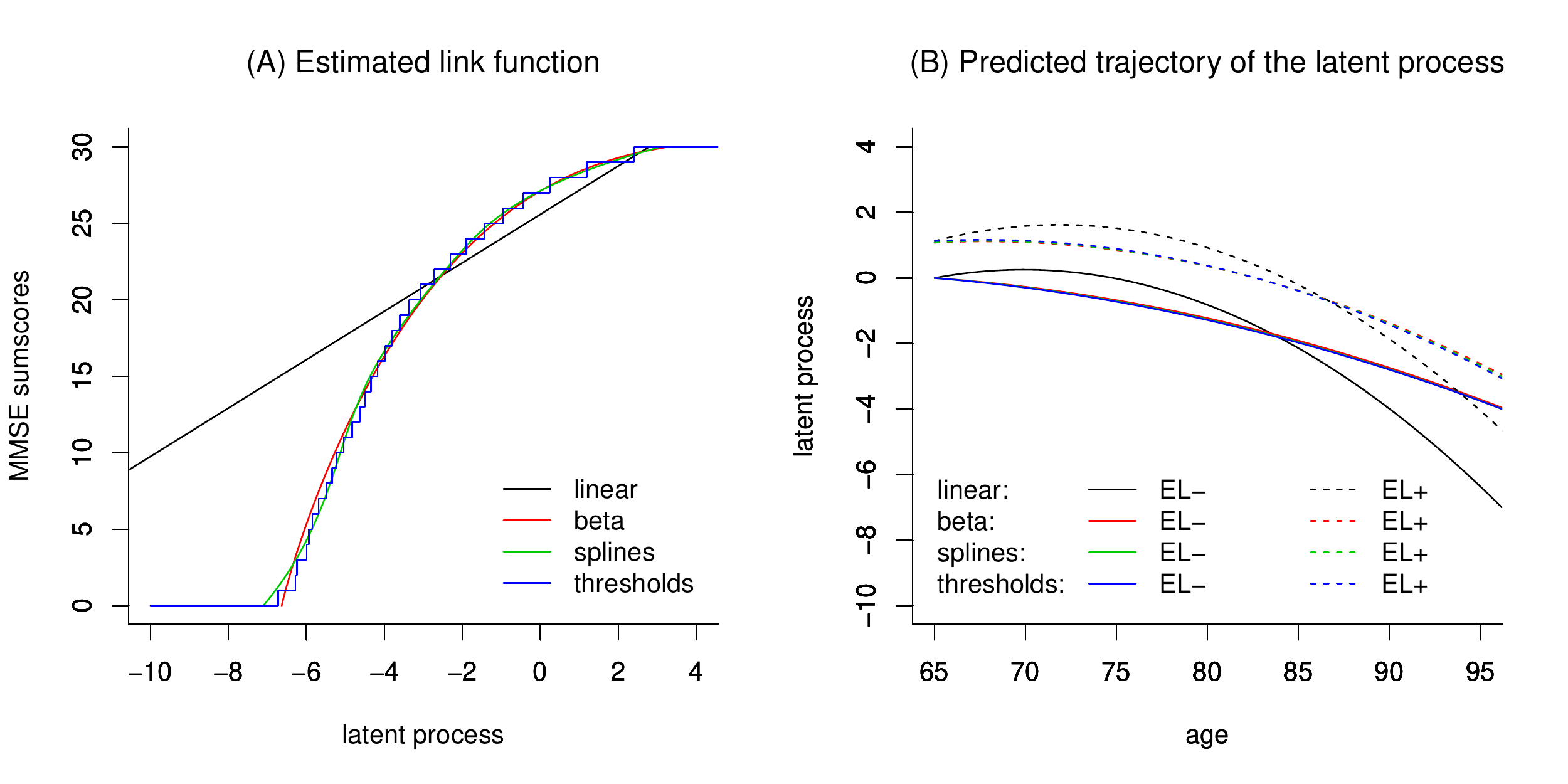}
\caption{(A) estimated inverse link functions between MMSE sumscores and the underlying latent process, and (B) predicted trajectories of the latent process of a woman according to educational level (with EL+ and EL- for respectively validated or non-validated primary school diploma) in latent process mixed models assuming either linear, Beta cdf, I-splines or threshold link functions (PAQUID sample, n=2,914); the trajectories for the latter three transformations are indistinguishable}
\label{histo_MMSE_CALC}
\end{figure}

For the Calculation subscore also, the standard linear mixed model again gave a clearly higher risk than the mixed model assuming a threshold link model ($D_{\UACV} {\rm (linear, Thresholds)}= 0.4523$, 95\% Tracking interval: $[0.4127,0.4919]$).

\begin{table}
\caption{Number of parameters (p), individual log-likelihood ($\Phi(\hat{\theta})$), corresponding naive normalized AIC (AIC), log-likelihood computed with respect to the counting measure ($\Psi(\hat{\theta})$), corresponding AIC (AIC$_d$), and UACV for latent process mixed models involving different transformations $H$ and applied on either the MMSE sumscore or its calculation subscore.}\label{table_UACV}
\begin{center}
\begin{tabular}{lcccccc}
\hline
Transformation $H$ & p & $\Phi(\hat{\theta})$ & AIC & $\Psi(\hat{\theta})$ & AIC$_d$ & UACV  \\
\hline
\multicolumn{7}{l}{\emph{MMSE}} \\
Linear & 16 & -8.7468 & 8.7523 & -8.5231 &8.5286&8.5361 \\
Beta cdf$^\dag$  & 18 & -7.7514 & 7.7576 & -7.7803& 7.7865& 7.7865\\
I-splines$^\ddag$ & 21 & -7.8249 & 7.8321 & -7.7857& 7.7929& 7.7935\\
Thresholds & 44 & -7.7473 & 7.7624 & -7.7473 & 7.7624 & 7.7665 \\
\hline
\multicolumn{7}{l}{\emph{Calculation}} \\
Linear & 16 & -6.0057 & 6.0111& -4.8143&4.8215 & 4.8198 \\
Thresholds & 19 & -4.3618  & 4.3683& -4.3618& 4.3683& 4.3692 \\
\hline
\end{tabular}
\end{center}
\vspace*{0.5cm}
$^\dag$ cdf for cumulative distribution function\\
$^\ddag$ Quadratic I-splines with 5 equidistant nodes located at 0, 7.5, 15, 22.5 and 30.  \\
\end{table}

\begin{table}
\caption{ Difference of AIC$_d$ ( $D_{\AIC_d}$ ), difference of UACV and its 95\% tracking interval between latent process mixed models involving different transformations $H_1$ and $H_2$, and applied on either the MMSE sumscore or its calculation subscore.}\label{table_diffUACV}
\begin{center}
\begin{tabular}{lrrr}
\hline
Transformations $H_1$ / $H_2$  & $D_{\AIC_d}$ & $D_{\UACV}$ & {95\% tracking interval} \\
\hline
\multicolumn{4}{l}{\emph{MMSE}} \\
linear / Beta cdf$^\dag$ &  0.7421 & 0.7495 & [0.6619 ; 0.8372] \\
linear / I-splines $^\ddag$ &  0.7357 & 0.7425 & [0.6526 ; 0.8325] \\
Beta cdf$^\dag$ / I-splines $^\ddag$ & -0.0064 & -0.0070 & [-0.0152 ;  0.0012]\\
I-splines $^\ddag$ / thresholds &  0.0306 & 0.0270 & [0.0166 ; 0.0374] \\
Beta cdf$^\dag$ / thresholds &  0.0241 & 0.0200 & [0.0097 ; 0.0303]\\
linear / thresholds &  0.7662 & 0.7696 & [0.6784 ; 0.8607] \\
\hline
\multicolumn{4}{l}{\emph{Calculation}} \\
linear / thresholds & 0.4515 & 0.4523 & [0.4127 ; 0.4919] \\
\hline
\end{tabular}
\end{center}
\vspace*{0.5cm}
$^\dag$ cdf for Cumulative Distribution Function\\
$^\ddag$ Quadratic I-splines with 5 equidistant nodes located at 0, 7.5, 15, 22.5 and 30.  \\
\end{table}


\section{Conclusion}
We have proposed a universal approximate formula for leave-one out crossvalidation: it is universal in the sense that it applies to any couple of estimating and assessment risks which can be correctly estimated from the observations. This is in principle restricted to parametric models but extends to smooth semi- or non-parametric ones through penalized likelihood. The approximate formula not only allows fast computation, but also allows deriving the asymptotic distribution. Estimating this distribution is important since the variability of UACV, as that of all the  criteria used for estimator choice, is large, even if the variability of a difference of UACV between two estimators is smaller.

In this paper, UACV has been applied to the issue of choice between estimators of the distribution of categorical data based on threshold models or on models based on a continuous approximation. It has been shown that the naive AIC can be misleading while a procedure called AIC$_d$ (which had not been theoretically validated) yields results very close to UACV, even if the latter is slightly better. Both quantities can be computed in the \texttt{lcmm} R package.

\bibliographystyle{chicago}

\bibliography{UACV}

\vspace{15mm}
\label{lastpage}
\end{document}